\documentclass[a4paper,12pt,centertags]{amsart}
\usepackage{pslatex}
\usepackage{dsfont}
\usepackage{hyperref}
\newtheorem{theorem}{Theorem}[section]
\newtheorem{lemma}[theorem]{Lemma}
\newtheorem{proposition}[theorem]{Proposition}

\theoremstyle{definition}
\newtheorem{definition}[theorem]{Definition}
\theoremstyle{remark}
\newtheorem{remark}[theorem]{Remark}
\numberwithin{equation}{section}

\DeclareMathOperator{\Div}{div}
\DeclareMathOperator{\supp}{Supp}
\newcommand{\cov}{\mathcal{O}}

\newcommand{\n}{\textbf{N}}
\newcommand{\Pro}{\mathds{P}}
\newcommand{\Es}{\mathds{E}}
\newcommand{\dint}{\int\!\!\!\!\int}
\newcommand{\intDt}{\int_0^T\!\!\!\!\int_D}
\newcommand{\e}{\varepsilon}
\hypersetup{pdftitle={Existence of martingale and stationary suitable weak solutions for a stochastic Navier-Stokes system},
pdfsubject={The existence of suitable weak solutions of 3D Navier-Stokes equations,
driven by a random body force, is proved. These solutions satisfy a local
balance of energy. Moreover it is proved also the existence of a 
statistically stationary solution.},
pdfauthor={F. Flandoli, M. Romito},
pdfkeywords={Navier-Stokes equations, suitable weak solutions martingale solutions, stationary solutions}}
\begin{document}
\title[Existence of martingale and stationary...]{Existence of martingale and stationary suitable weak 
       solutions for a stochastic Navier-Stokes system}
\author{Marco Romito}
\address{Dipartimento di Matematica \emph{U. Dini}, Universit\`a di
Firenze, Viale Morgagni 67/a, 50134 Firenze, Italia}
\email{romito@math.unifi.it}
\subjclass{Primary 76D05; Secondary 35R60}
\date{January 11, 2001}
\keywords{Navier-Stokes equations, suitable weak solutions martingale
solutions, stationary solutions}
\begin{abstract}
The existence of suitable weak solutions of 3D Navier-Stokes equations,
driven by a random body force, is proved. These solutions satisfy a local
balance of energy. Moreover it is proved also the existence of a 
statistically stationary solution.
\end{abstract}
\maketitle

\section{Introduction}
The suitable weak solutions of three dimensional Navier-Stokes system are 
weak solutions which satisfy a local energy inequality. The local energy 
inequality can be seen as a mathematical counterpart of a local balance 
of the variation of the energy of the fluid governed by the Navier-Stokes 
equations. This additional property in general is not verified by solutions 
to Navier-Stokes equations (e.g. solutions obtained by Galerkin approximations),
but, from a physical viewpoint, it is natural to choose solutions which have 
more regularity properties and so a more precise physical meaning.

The concept of suitable weak solution was introduced firstly by Caffarelli, 
Kohn and Nirenberg in \cite{CKN}, with the aim to study the partial regularity 
of solutions of Navier-Stokes system. With this approach the local energy 
inequality is turned out to be a crucial tool and it has brought Caffarelli, 
Kohn and Nirenberg to show in \cite{CKN} the best known result in this setting.

The suitable weak solutions were already defined in the papers of Scheffer 
\cite{Sc1}, \cite{Sc2}, \cite{Sc3}, even if in a rather implicit way. 
Scheffer gives an interpretation of the local energy inequality in terms 
of the existence of an external force $f$ acting \textsl{against} the flow, 
in other words
$$
f(t,x)\cdot u(t,x)\le 0.
$$

In this paper we will show the existence of martingale suitable weak 
solutions of three dimensional Navier-Stokes system. These are solutions 
of the Navier-Stokes equations  driven by a random body force (a white 
noise). Different interpretations can be given to these terms. A random 
force can represent all those phenomena that are usually neglected where 
the system is stable. A random solution of the Navier-Stokes equations 
can take into account that flows with a large Reynolds number are chaotic 
and sensitive to microscopic perturbations. Sources of perturbations can 
be the roughness of the boundary, or the effects of the environment 
external to the system, such as acoustic waves, etc. 

We will prove also the existence of stationary suitable weak 
solutions, where stationary has to be understood in a 
statistical sense. These solutions may describe a fluid in a turbulent 
regime. The partial regularity of these solutions will be studied in 
other papers (see \cite{FlRo1}, \cite{FlRo2} and \cite{Rom}).

Many authors have proved the existence of deterministic suitable weak 
solutions, see for example Caffarelli, Kohn and Nirenberg \cite {CKN}, 
or P.L. Lions \cite{Lio}, Beirao da Veiga \cite{Bei}, Lemarie-Rieusset 
\cite{Lem}. None of the methods used in these papers is well suited to 
prove the existence of martingale solutions. The existence of suitable 
weak solutions in a stochastic setting seems to be new.

\subsection{Notations}
Let $D$ be an open bounded subset of ${\bf R}^3$ having a smooth boundary 
and for each $T>0$ set $D_T=(0,T)\times D$. Define the space
$$
H=\{u:D\to{\bf R}^3\,|\,u\in\left(L^2(D)\right)^3,\,{\mathop{\rm div}} u=0,\,u\cdot n|_{\partial D}=0\,\},
$$
where $n$ is the outer normal to $\partial D$ (see for example Temam 
\cite{Te1}), and the space
$$
V=\{u\in \left(H^1(D)\right)^3\,|\,{\mathop{\rm div}} u=0,\,u|_{\partial D}=0\,\}.
$$
The $L^2$-norm of elements of $H$ will be denoted by $|\cdot|$ and 
the $H^1$-norm of elements of $V$ will be denoted by $\|\cdot\|$. 
If the space $H$ is identified with its dual and the space $H'$ is 
identified with a subspace of $V'$, then
$$
V\subset H\subset V'.
$$
The operator $A:D(A)\subset H\to H$ is defined as 
$Au=-{\mathcal P}\triangle u$, where ${\mathcal P}$ is the orthogonal 
projection from $L^2(D)^3$ onto $H$ and $D(A)=H^2(D)^3\cap V$. The 
operator $A$ is positive self adjoint with compact resolvent. The 
eigenvalues of $A$ are denoted by $0<\lambda_1\le\lambda_2\le\ldots$ 
and $e_1$, $e_2$, \dots is a corresponding complete orthonormal 
system of eigenvectors.

Moreover, if the dual space of $D(A)$ is denoted by $D(A^{-1})$, then 
identifications as above give the dense continuous inclusions
$$
D(A)\subset V\subset H\subset V'\subset D(A^{-1}).
$$

The fractional powers $A^\alpha$ of $A$, $\alpha\ge 0$, are simply 
defined by
$$
A^\alpha x=\sum_{i=1}^\infty\lambda_i^\alpha\langle x,e_i\rangle e_i
$$
with domain
$$
D(A^\alpha)=\{\, x\in H\,|\, \|x\|_{D(A^\alpha)}<\infty\,\}
$$
where
$$
\|x\|_{D(A^\alpha)}^2=\sum_{i=1}^\infty\lambda_i^{2\alpha}\langle x,e_i\rangle^2=|A^\alpha x|^2.
$$
The space $D(A^\alpha)$ is an Hilbert space with the inner product
$$
\langle x,y\rangle_{D(A^\alpha)}=\langle A^\alpha x,A^\alpha y\rangle,\qquad x,y\in D(A^\alpha).
$$
Notice that $D(A^\alpha)\subset H^{2\alpha}(D)$.

Since $V$ coincides with $D(A^{1/2})$ (see Temam \cite{Te2} Section 2.2, 
or Temam \cite{Te3}, Ch. III, Section 2.1), the space $V$ can be endowed 
with the norm $\|u\|=|A^{1/2}u|$. The Poincar\`e inequality gives
\begin{equation*}
||u||^{2}\geq \lambda _{1}|u|^{2}.
\end{equation*}

Moreover we will consider the Sobolev spaces $W^{s,p}(0,T;H)$ endowed 
with the norm
$$
\|f\|_{W^{s,p}}^p=\int_0^T|f|^p\,dt+\int_0^T\!\!\!\int_0^T\frac{|f(t)-f(s)|^p}{|t-s|^{1+sp}}\,ds\,dt
$$

Finally we define the bilinear operator $B:V\times V\to V'$ as
$$
\langle B(u,v),w\rangle=\int_D w\cdot (u\cdot\nabla)v\,dx,\qquad w\in V.
$$
The operator $B$ can be extended in many different spaces (see for 
example Temam \cite{Te2}).

In the sequel we will largely use the following {\sl local} Sobolev 
inequality. Let $u\in H^1(B_r)$, then
\begin{equation}\label{sob}
\int_{B_r}|u|^q\le C\left(\int_{B_r}|\nabla u|^2\right)^a\left(\int_{B_r}|u|^2\right)^{\frac{q}2-a}+\frac{C}{r^{2a}}\left(\int_{B_r}|u|^2\right)^{\frac{q}2},
\end{equation}
where $q\in[2,6]$ and $a=\frac34(q-2)$.


\section{Definitions and main results}

We consider a viscous incompressible homogeneous Newtonian fluid in the 
bounded open domain $D\subset{\bf R}^3$, described by the Navier-Stokes 
equations
\begin{equation}\label{NS}
\begin{cases}
\partial_t u+(u\cdot\nabla)u+\nabla P=\nu\triangle u+f+\partial_tg 
   &\qquad\text{in }D_T\\
\Div u=0
   &\qquad\text{in }D_T\\
u=0
   &\qquad\text{on }[0,T]\times\partial D\\
u(0)=u_0
\end{cases}
\end{equation}
where $u$ is the velocity field, $P$ is the pressure field and $\nu$ is 
the kinematic viscosity. For simplicity we will take $\nu=1$, since its 
value is not relevant in the present setting. The term $\partial_tg$ 
represents a rapidly fluctuating force and in this paper it will be 
model by a white noise. In order to handle this term we introduce 
the new variables
$$
v=u-z,\qquad \pi=P-Q,
$$
where the pair $(z,Q)$ solves the following Stokes equation
\begin{equation}\label{lNS}
\begin{cases}
\partial_t z+\nabla Q=\triangle z+f+\partial_t g,\\
\Div z=0,\\
z=0&\qquad\text{on $\partial D$},\\
z(0)=0,
\end{cases}
\end{equation}
in $[0,T]\times D$. Then the new variables $(v,\pi)$ solve the 
following equation
\begin{equation}\label{mNS}
\begin{cases}
\partial_t v+\left((v+z)\cdot\nabla\right)(v+z)+\nabla\pi=\triangle v,\\
\Div v=0,\\
v=0&\qquad\text{on $\partial D$},\\
v(0)=u_0,
\end{cases}
\end{equation}
where the white noise term has disappeared. We will ask for a path-wise 
local energy inequality to the pair $(v,\pi)$, as we shall see in the 
sequel.


\subsection{Assumptions on the data} 

We will model the fluctuation part of the body force $\partial_t g$
as a noise white in time $\partial_t B$, so that $B$ is a Brownian motion.
We will assume throughout the paper the following assumptions
\begin{equation}\label{As}
\begin{aligned}
&u_0\in H,\\
&f\in L^2(0,T;H),\\
&B\text{ is a Brownian motion with trajectories in }D(A^\delta)
\end{aligned}\tag{As}
\end{equation}
for a small $\delta>0$. It is possible to see that this set of assumptions 
implies that the trajectories of the solution $z$ of equations \eqref{lNS}
have the following regularity properties
\begin{equation}\label{zreg}
z\in L^\infty(0,T;H)\cup L^2(0,T;V)\cap L^\infty(0,T;L^4(D))\qquad\Pro-\text{a.s.}
\end{equation}
(see Flandoli \cite{Fla}).

\begin{remark}
Another way to give assumptions \eqref{As} regarding the Brownian motion
$B$ is to consider the covariance operator $\cov$, which is a positive bounded
self-adjoint operator. We suppose that $\cov$ maps $H$ into $D(A^\delta)$. So
a sufficient condition that ensures \eqref{As} is that the operator
$A^\delta\cov A^\delta$ has a bounded extension to $H$ which is of trace class
(see \cite{DaPZa} for more).
\end{remark}

\begin{remark}
Even if we are mainly interested in interpreting the fluctuation $\partial_t g$
as a white-noise, in view of stationary solutions we will consider also
deterministic solutions. In this case we will assume
\begin{equation}\label{Ad}
\begin{aligned}
&u_0\in H,\\
&f\in L^2(0,T;H),\\
&g\in C^{\frac12-\e}([0,T];D(A^\delta)),\qquad g(0)=0,
\end{aligned}\tag{Ad}
\end{equation}
for $\delta>\e>0$, so that again property \eqref{zreg} holds (we refer again 
to Flandoli \cite{Fla}). Note that it is possible to choose the function $g$
with different regularity properties, combining in different ways the
differentiability with respect to time and the differentiability with respect
to space.
\end{remark}


\subsection{Martingale suitable weak solutions}

We start with the definition of martingale suitable weak solutions and
we give the main theorem about their existence.

Before doing this, we define the suitable solutions in a deterministic
setting, so that the derivative of $g$ with respect to time has to be
understood in the sense of distributions.

\begin{definition}\label{detsws}
Let $T\in(0,\infty]$. A suitable weak solution to Navier-Stokes equations
is a pair $(u,P)$ such that if $v=u-z$ and $\pi=P-Q$, where $(z,Q)$ is the
solution of equation \eqref{lNS}, then
\begin{enumerate}
\item $v$ is weakly continuous with respect to time,
\item $v\in L^\infty(0,T;H)\cap L^2(0,T;V)$ and $\pi\in L^{5/3}_{\text{loc}}(D_T)$,
\item $(v,\pi)$ satisfies equation \eqref{mNS} in the sense of distributions on $D_T$,
\item for all $t\le T$ and almost all $s<t$,
$$
|v(t)|^2+2\int^t_s\|v\|^2\,dr
\le |v(s)|^2+\int^t_s\int_D z\cdot\left((v+z)\cdot\nabla\right)v\,dx\,dr
$$
\item [{\it (v)}] for any $\phi\in C^\infty_c(D_T)$, $\phi\ge 0$,
\begin{align}\label{lei}
2\intDt |\nabla v|^2\phi
&\le     \intDt |v|^2(\partial_t\phi+\triangle\phi)\notag
        +2\intDt \pi(v\cdot\nabla\phi)\\
&\quad  +\intDt (|v|^2+2v\cdot z)\left((v+z)\cdot\nabla\phi\right)\notag\\
&\quad  +2\intDt \phi z\cdot\left((v+z)\cdot\nabla\right)v\notag
\end{align}
\end{enumerate}
\end{definition}

A martingale suitable weak solution for the Navier-Stokes equations will 
be the solution of a stochastic differential equation driven by an additive
noise such that its trajectories are suitable weak solutions in the sense 
of the definition above. More precisely:

\begin{definition}\label{marsws}
A martingale suitable weak solution is a process $(u,P)$ defined on a
stochastic basis
$$
(\Omega,\mathcal{F},(\mathcal{F})_{t\ge 0},\Pro,(B_t)_{t\ge0}),
$$
where $B$ is a Brownian motion adapted to the filtration with values in
$D(A^\delta)$, such that
$$
\omega\in\Omega\to(u(\omega),P(\omega))\in L^2(0,T;H)\times L^{5/3}_\text{loc}(D_T)
$$
is a measurable mapping and such that there exists a set $\Omega_0\subset\Omega$
of full probability such that the pair $(u(\cdot,\omega),P(\cdot,\omega))$ is a
suitable weak solution in the sense of Definition \ref{detsws}, with
respect to the body force $f+\partial_t B_t(\omega)$, for all
$\omega\in\Omega_0$.
\end{definition}

We want to explain the meaning of the last part of this definition. Since
$B$ is a Brownian motion, under assumption \eqref{As} it has $\Pro$-a.s.
trajectories in $C^{\frac12-\e}([0,T];D(A^\beta))$, for
$0<\e<\beta\le\delta$. So for every given $\omega\in\Omega_0$, we have
that $B(\omega)\in C^{\frac12-\e}([0,T];D(A^\beta))$, the solution
$(z(\omega),Q(\omega))$ enjoys the regularity stated in \eqref{zreg} and
the pair $(u(\omega),P(\omega))$ satisfies all the conditions of
Definition \ref{detsws} with respect to these functions.

We can give now the main existence theorem for martingale suitable weak
solutions.

\begin{theorem}\label{thmar}
Assume \eqref{As}. There exists a martingale suitable weak solution
(in the sense of Definition \ref{marsws} above) with initial data $u_0$.
Moreover
\begin{equation}\label{energyinequality}
{\Es}|u(t)|^2_H+{\Es}\int_s^t\|u\|_V^2\,dr
\le {\Es}|u(s)|^2_H 
    +\sigma(t-s)
    +{\Es}\int^t_s\|f\|_{V'}\,ds,
\end{equation}
and
\begin{align}\label{supenergy}
\Es\bigl[\sup_{(s,t)}|u_N(r)|^2_H\bigr]+\Es\int_s^t\|u_N\|_V^2\,dr
&\le     2{\Es}|u_0|^2_H
       + 2\int^t_0\|f\|_{V'}^2\,ds\notag\\
&\quad + 2\sigma(1+\sigma C_1^2) (t-s).
\end{align}
where $\sigma$ is the variance of $B$ and $C_1$ is a universal constant.
\end{theorem}

\begin{remark}
It can be noticed that, as in \cite{CKN}, the complete local energy inequality 
\begin{align*}
&\int_D|u(t)|^2\phi+2\int^t_0\int_D |\nabla v|^2\phi\le\\
&\qquad\le    \int^t_0\int_D |v|^2(\partial_t\phi+\triangle\phi)
             +2\int^t_0\int_D \pi(v\cdot\nabla\phi)\\
&\qquad\quad +\int^t_0\int_D (|v|^2+2v\cdot z)\left((v+z)\cdot\nabla\phi\right)
             +2\int^t_0\int_D \phi z\cdot\left((v+z)\cdot\nabla\right)v
\end{align*}
can be recovered using a {\sl cut-off} function $\chi$:
$$
0\le\chi\le1,\qquad \chi=0\mbox{ for }s\le 0,\qquad \chi=1\mbox{ for }s\ge1;
$$
for each $t$ we use $\phi_\e(x,s)=\phi(x,s)\chi(\frac{t-s}{\e})$
as a test function and, as $\e\to 0$, we obtain the full local
energy inequality.
\end{remark}

\begin{remark}
The definition of suitable weak solution we have given seems to depend
on the solution $z$ of the linear problem. This is not true, the
definition given above has been introduced only to deal with the term
$\partial_t g$. Indeed, it is possible to show the following result, which
will be proved in Section \ref{secnozeta}.

\begin{theorem}\label{nozeta}
The property of being a suitable weak solution for a pair $(u,P)$ does
not depend on the solution $(z,Q)$ chosen for the linear problem.
\end{theorem}

The previous theorem tells us that, when $g\equiv0$, there is no difference
between the suitable weak solutions in the sense of Caffarelli, Kohn and
Nirenberg \cite{CKN} and ours.
\end{remark}


\subsection{Stationary solutions}

The approach we follow here concerning the framework of the path space
and the introduction of stationary solutions is due to Sell \cite{Sel}
(see also \cite{FlSc}) and gives a solution to the problem of studying
the asymptotic behaviour of dynamics when the dynamic itself cannot be
well defined, as for Navier-Stokes equations.

A stationary solution is a measure on the space of all trajectories 
$(u,W)$ that are solutions to Navier-Stokes equations, which is invariant
for the time-shift. In this setting we will not consider the pressure term 
$P$ explicitly, since we are mainly interested in the statistical properties
of the velocity. In fact in \cite{FlRo2} a regularity criterion will be 
proved which involves only the gradient of the velocity.

In order to have an equation whose deterministic part is autonomous, we 
will suppose that the deterministic forcing term $f\in L^2(D)$ is 
independent of time. The time-shift will act on the increments of the 
Brownian motion in order to preserve the stationarity of its increments.

Let $C_0([0,+\infty),H)$ be the set of all continuous functions which
take value $0$ in $t=0$ and let ${\mathcal S}$ be the subset of
$L^2_{\rm loc}(0,+\infty;H)\times C_0([0,+\infty),H)$ of all suitable
weak solutions in $(0,\infty)\times D$, that is the set of all pairs
$(u,W)$, where $W\in C^{1/2-\e}([0,T];D(A^\beta))$ for 
$0<\e<\beta\le\delta$, and $u$ is a suitable weak solution in the
sense of Definition \ref{detsws} for all $T>0$ under the body force
$f+\partial_t W$. In this setting the pressure $P$ is treated as an
auxiliary scalar field. We will see that this set is not empty. Let 
us define a metric on ${\mathcal S}$. Let
\begin{align*}
d_1(u^1,u^2)
&=\sum_{n=1}^\infty 2^{-n}\Bigl(1\wedge\int_0^n|v^1-v^2|^2\,dt\Bigr)^{\frac12},\\
d_2(W^1,W^2)
&=\sum_{n=1}^\infty 2^{-n}\bigl(1\wedge\sup_{(0,n)}|W^1-W^2|\bigr),
\end{align*}
and the metric on ${\mathcal S}$ is defined as
$$
d\left((u^1,W^1),(u^2,W^2)\right)=d_1(u^1,u^2)+d_2(W^1,W^2).
$$

Let $C_b({\mathcal S})$ be the space of all bounded real continuous 
functions on ${\mathcal S}$ with the uniform topology, let $\mathcal B$ 
be the Borel $\sigma$-algebra of $({\mathcal S},d)$ and 
$M_1({\mathcal S})$ be the set of all probability measures on 
$({\mathcal S},{\mathcal B})$.

Let $\tau_t:{\mathcal S}\to {\mathcal S}$, ($t\ge 0$) be the time 
shift on ${\mathcal S}$, defined as
$$
\tau_t(u,W)(s)=(u(s+t),W(t+s)-W(t)).
$$
Notice that the map $(t,u,W)\to\tau_t(u,W)$ is continuous from 
$[0,\infty)\times {\mathcal S}$ to ${\mathcal S}$. We denote again 
by $\tau_t$ the induced mapping on $C_b({\mathcal S})$, defined as
$$
\tau_t\phi(u)=\phi(\tau_t u)
$$
and by $\tau_t\mu$ the image measure of any $\mu\in M_1({\mathcal S})$ 
under $\tau_t$, in the sense that
$$
\langle\tau_t\mu,\phi\rangle=\langle\mu,\tau_t\phi\rangle
$$
for each $\phi\in C_b({\mathcal S})$.

\begin{definition}
A probability measure $\mu\in M_1({\mathcal S})$ is time-stationary 
if $\tau_t\mu=\mu$ for all $t\ge 0$. A probability measure $\mu$ has 
finite mean dissipation rate if
$$
\int_{\mathcal S}\left[\int_0^T\int_D|\nabla u|^2\,dx\,dt\right]\mu(du)<\infty
$$
for all $T>0$.
\end{definition}

\begin{remark}
The property of having a finite dissipation rate is exactly the one we 
will need to apply the regularity criterion presented in \cite{FlRo2}. 
Notice that this property does not depend on the presence of the noise, 
since there exist stationary solutions with finite dissipation rate also 
for the deterministic solution (see \cite{FlRo1}).
\end{remark}

\begin{theorem}\label{thstat}
Let $f\in L^2(D)$ be independent of time. There exists a time stationary
probability measure $\mu\in M_1(\mathcal{S})$ with finite mean dissipation
rate. Moreover, there exists a constant $C_\mu>0$ such that for all
$t\ge s\ge0$,
\begin{equation}\label{lineare}
\int_{\mathcal S}\left[\int_0^T\int_D|\nabla u|^2\,dx\,dt\right]\mu(du)
=C_\mu(t-s).
\end{equation}
Finally, the image measure of $\mu$ under the projection onto the second
component is a Wiener measure whose covariance operator maps $H$ in
$D(A^\delta)$, for a small $\delta>0$.
\end{theorem}

The last claim of the theorem says poorly that the standard process on
$\mathcal{S}$ having law $\mu$ is a martingale suitable weak solution
driven by a Brownian motion satisfying assumption \eqref{As}.


\section{The proof of Theorem \ref{nozeta}}\label{secnozeta}

Let $(u,P)$ be a suitable weak solution in the sense of Definition
\ref{detsws}, that is, with respect to the solution $(z,Q)$ of
problem \eqref{lNS}. Let $(z_1,Q_1)$ be the solution of the equation
$$
\begin{cases}
\partial_t z_1-\triangle z_1+\nabla Q_1=f_1,\\
\Div z_1=0,
\end{cases}
$$
with initial condition $z_1(0)=z_0$, and set
$$
w=z-z_1\qquad\text{and}\qquad R=Q-Q_1.
$$
Then $v_1=v+w$, $\pi_1=R+\pi$ and $v+z=v_1+z_1=u$. We show that $(u,P)$
is a suitable weak solution with respect to $(z_1,Q_1)$. In order to
do this, we have only to show that $v_1$ satisfies the local energy
inequality \eqref{lei1} (which has an additional term which takes into
account the term $f-f_1$).

The function $w$ is the solution of
$$
\begin{cases}
\partial_t w-\triangle w+\nabla R=f_2,\\
\Div w=0,\\
w(0)=-z_0,
\end{cases}
$$
where $f_2=f-f_1$. It is easy to see, by mollification, that for any
$\phi\in C^\infty_c(D_T )$, $\phi\ge 0$,
$$
2\intDt |\nabla w|^2\phi
=\intDt |w|^2(\partial_t\phi+\triangle\phi)
 +2\intDt  R(w\cdot\nabla\phi)
 +2\intDt\phi f_2\cdot w.
$$

\begin{lemma}
With the notations above, for any $\phi\in C^\infty_c(D_T )$,
\begin{align*}
&4\intDt \phi\nabla v\cdot\nabla w=\\
&\qquad=      2\intDt v\cdot w(\partial_t\phi+\triangle\phi)
             +2\intDt R(v\cdot\nabla\phi)
             +2\intDt\pi(w\cdot\nabla\phi)\\
&\qquad\quad +2\intDt\phi f_2\cdot v
             +\intDt (2w\cdot z-|w|^2)\left((v+z)\cdot\nabla\phi\right)\\
&\qquad\quad +2\intDt \phi z\cdot\left((v+z)\cdot\nabla\right)w
             -2\intDt \phi w\cdot\left((v+z)\cdot\nabla\right)v_1
\end{align*}
\end{lemma}
\begin{proof} Let $\phi\in C^\infty_c(D_T )$; by mollification
in a neighbourhood $U$ of $\supp\phi$ we obtain $(w_\e,R_\e)$ 
such that $w_\e\to w$ in $L^\infty(L^2)$, 
$\nabla w_\e\to\nabla w$ in $L^2$ and $R_\e\to R$ in 
$L^{5/3}$ in $U$.

Since $(v,\pi)$ is a weak solution of \eqref{mNS}, we use 
$\phi w_\e$ as a test function to have
\begin{multline}\label{rel1}
   \intDt v\cdot w_\e\partial_t\phi
+\intDt \phi v\cdot\partial_t w_\e
+\intDt (v+z)\cdot\left[(v+z)\cdot\nabla\right](\phi w_\e)+\\
+\intDt \pi w_\e\cdot\nabla\phi
 = \intDt\nabla v\cdot\nabla(\phi w_\e)
\end{multline}
Moreover we have in $\supp\phi$
$$
\begin{cases}
\partial_t w_\e-\triangle w_\e+\nabla R_\e=f_2,\\
\Div w_\e=0
\end{cases}
$$
and multiplying by $\phi v$ and integrating by parts gives
\begin{equation}\label{rel2}
 \intDt\phi v\cdot\partial_t w_\e
+\intDt \nabla w_\e\cdot\nabla(\phi v)
=\intDt R_\e(v\cdot\nabla\phi)
+\intDt \phi f_2\cdot v.
\end{equation}
Then we subtract \eqref{rel2} from \eqref{rel1} and we use the following
facts (they can be easily obtained by integration by parts)
\begin{align*}
\intDt \nabla v\cdot\nabla(\phi w_\e)
  +\nabla w_\e\cdot\nabla(\phi v)
&=\intDt 2\phi\nabla v\cdot\nabla w_\e
 -v\cdot w_\e\triangle\phi,\\
\intDt \phi v\big((v+z)\cdot\nabla\big)w_\e
&=     -\intDt \phi w_\e\big((v+z)\cdot\nabla\big)v+\\
&\quad -\intDt (v\cdot w_\e)\big((v+z)\cdot\nabla\phi\big)\\
2\intDt \phi w_\e\big((v+z)\cdot\nabla\big)w_\e
&=-\intDt |w_\e|^2(v+z)\cdot\nabla\phi,
\end{align*}
so that we finally have
\begin{align*}
&4\intDt \phi\nabla v\cdot\nabla w_\e=\\
&\qquad  2\intDt v\cdot w_\e(\partial_t\phi+\triangle\phi)
        +2\intDt R_\e\,v\cdot\nabla\phi
        +2\intDt\pi\,w_\e\cdot\nabla\phi\\
&\qquad +2\intDt\phi v\cdot f_2
        +\intDt(2z\cdot w_\e-|w_\e|^2)\left[(v+z)\cdot\nabla\phi\right]\\
&\qquad +2\intDt\phi z\left[(v+z)\cdot\nabla\right]w_\e
        -2\intDt\phi w_\e\left[(v+z)\cdot\nabla\right](v+w_\e)
\end{align*}
and, as $\e$ goes to $0$, the conclusion follows.
\end{proof}

With the help of the above lemma, we can conclude the proof of Theorem
\ref{nozeta}. Let $\phi\in C^\infty_c(D_T )$, with $\phi\ge 0$, and 
$t\in(0,T]$. By definition
$$
2\intDt \phi|\nabla v_1|^2
= 2\intDt \phi|\nabla v|^2
 + 2\intDt \phi|\nabla w|^2\\
 + 4\intDt \phi\nabla v\cdot\nabla w.
$$
Since $(v,\pi)$ satisfies the energy inequality and using the energy 
equality for $(w,R)$ and the previous lemma, we have
\begin{align*}
&2\intDt \phi|\nabla v_1|^2\le\\
&\qquad\le    \intDt (|v|^2+2v\cdot w+|w|^2)(\partial_t\phi+\triangle\phi)
             +2\intDt (\pi+R)(v\cdot\nabla\phi)\\
&\qquad\quad +\intDt (|v|^2+2v\cdot z+2w\cdot z-|w|^2)\left((v+z)\cdot\nabla\phi\right)\\
&\qquad\quad +2\intDt \phi z\cdot\left((v+z)\cdot\nabla\right)(v+w)
             +2\intDt (R+\pi)(w\cdot\nabla\phi)\\
&\qquad\quad +2\intDt \phi f_2\cdot(v+w)
             -2\intDt \phi w\cdot\left((v+z)\cdot\nabla\right)v_1.
\end{align*}
Now we use the fact that $v_1=v+w$, $\pi_1=\pi+R$, $v+z=v_1+z_1$ and that
$$
|v|^2+2v\cdot z+2w\cdot z-|w|^2=|v_1|^2+2v_1\cdot z_1
$$
to obtain the local energy inequality for $v_1$
\begin{align}\label{lei1}
2\intDt \phi|\nabla v_1|^2
&\le    \intDt |v_1|^2(\partial_t\phi+\triangle\phi)
       +2\intDt \pi_1(v_1\cdot\nabla\phi)\notag\\
&\quad +\intDt (|v_1|^2+2v_1\cdot z_1)\left((v_1+z_1)\cdot\nabla\phi\right)\\
&\quad +2\intDt \phi z_1\cdot\left((v_1+z_1)\cdot\nabla\right)v_1
       +2\intDt\phi v_1\cdot f_2.\notag
\end{align}


\section{Proofs of the existence theorems}

In this section we will prove Theorem \ref{thmar}, on existence for
martingale suitable weak solutions, and Theorem \ref{thstat}, about
the existence of stationary solutions. Prior to do this, we show
a path-wise existence result, which will be the basis of the proofs
of the two main theorems.


\subsection{Path-wise existence}

In this section it will be proved the existence of suitable weak solutions
as defined in Definition \ref{detsws}. In other words we will show
the following theorem

\begin{theorem}\label{thdetsws}
Assume \eqref{Ad}. There exists a suitable weak solution in the sense
of Definition \ref{detsws}.
\end{theorem}

The proof of this theorem is given in three steps. In the first step we
solve a linearised version of the equation, whose higher regularity
will let us prove the local energy inequality for such solutions. The 
second step will consist in the application of the Banach fixed point
theorem to get the solution of an approximated nonlinear problem.
In the third step we will find, in the limit of the approximation,
a solution as requested by Theorem \ref{thdetsws}.

We start with the first step.

\begin{lemma}\label{exlemma}
Let $u_0\in V$, $w\in L^\infty\left(D_T \right)\cap L^2(0,T;V)$
and $\xi\in L^2(0,T;V)\cap L^\infty(0,T;H)$. Then there exists
a unique solution $(u,p)$ of the problem
$$
\begin{cases}
\partial_t u-\triangle u+\nabla p+(w\cdot\nabla)(u+\xi)=0,\\ 
\Div u=0,\\ 
u=0&\qquad\text{on }(0,T)\times\partial D,\\
u(0)=u_0,
\end{cases}
$$
with $u\in C([0,T];H^1_0(D))\cap L^2(0,T;H^2(D))$ and 
$p\in L^2(0,T;H^1(D))$.

Moreover we have  $\partial_t u\in L^2(D_T)$ and for
any $\phi\in C^\infty_c(D_T )$, $\phi\ge 0$,
\begin{align*}
2\intDt |\nabla u|^2\phi
&=     \intDt |u|^2(\partial_t\phi+\triangle\phi)
      +2\intDt  p(u\cdot\nabla\phi)\\
&\quad+2\intDt (|u|^2+2u\cdot\xi)(w\cdot\nabla\phi)
      +2\intDt  \phi\xi\cdot(w\cdot\nabla)u
\end{align*}
\end{lemma}
\begin{proof}
Let ${\mathcal B}:L^2(0,T;V)\to L^2(0,T;L^2(D))$ be the operator $$
{\mathcal B}u=(w\cdot\nabla)u $$ It is easy to see that $({\mathcal
B}u,v)_{L^2(D)}=-({\mathcal B}v,u)_{L^2(D)}$ for $u$, $v\in V$ and so
it is easy to deduce that $({\mathcal B}u,u)_{L^2(D)}=0$.

We look for a solution $u\in L^2(0,T;V)$ of the problem
$$
\begin{cases}
\frac{d}{dt}(u,v)_H+(u,v)_V+({\mathcal B}(u+\xi),v)_{L^2(D)}=0
   &\qquad\text{for each $v\in V$},\\
u(0)=u_0.
\end{cases}
$$
We will see that $u\in V\cap H^2$, then $Au\in H$ and so
$u'=-Au-{\mathcal PB}(u+\xi)$ in H and $u'\in L^2(0,T;H)$, i.e. $u$ is
equal a.e. to a continuous function from $[0,T]$ to $H$ (see
Temam \cite{Te1}, Lemma 3.1.1) and the initial condition makes sense.

We prove existence by means of the Galerkin method. Let
$v_1,\ldots,v_m,\ldots$ be a basis as above. We define for each $m\ge 1$ the
approximate solution $u_m$ of the problem as follows
$$
u_m=\sum_{i=1}^m u_m^i(t)v_i
$$
and
$$
\begin{cases}
(u'_m,v_j)_H+(u_m,v_j)_V+({\mathcal B}(u_m+\xi),v_j)_{L^2(D)}=0
   &\qquad j=1,\ldots,m\\
u_m(0)=u_m^0,
\end{cases}
$$
where $u_m^0$ is the orthogonal projection in $H$ of $u_0$ on the linear 
space spanned by $v_1,\ldots,v_m$. This finite-dimensional linear system 
has a unique solution.

First we obtain an estimate of $u_m$ in $L^\infty(0,T;H)\cap
L^2(0,T;V)$. Multiply each equation respectively by $u_m^i$ and
sum to have
$$
\frac{d}{dt}|u_m|_H^2+2||u_m||_V^2+2\left({\mathcal B}(u_m+\xi),u_m\right)_{L^2(D)}=0.
$$
Since
\begin{align*}
2\left({\mathcal B}(u_m+\xi),u_m\right)_{L^2(D)}
& = -2\left({\mathcal B}u_m,\xi\right)_{L^2(D)}\\
&\le \|u_m\|_V^2+\int_D |\xi|^2|w|^2\\
&\le \|u_m\|_V^2+\|w\|_{L^\infty(D)}^2\|\xi\|_{L^2(D)}^2,
\end{align*}
we have
$$
\frac{d}{dt}|u_m|_H^2+||u_m||_V^2\le \|w\|_{L^\infty(D)}^2\|\xi\|_{L^2(D)}^2
$$
and consequently, integrating in time,
$$
\sup_{(0,T)}|u_m|^2_H\le
|u_m^0|_H^2+T\|w\|^2_{L^\infty(D_T)}\|\xi\|^2_{L^\infty(L^2(D))}
$$
and
$$
\int^T_0\|u_m\|_V^2\le
|u_m^0|_H^2+T\|w\|^2_{L^\infty(D_T)}\|\xi\|^2_{L^\infty(L^2(D))}.
$$

Then we obtain an estimate of $u_m$ in $L^\infty(0,T;V)\cap
L^2(0,T;H^2( D ))$ and of $u'_m$ in $L^2(0,T;H)$. Multiply each
equation by ${u_m^i}'$ and sum to obtain
$$
|u'_m|_H^2+\frac12\frac{d}{dt}||u_m||^2_V+({\mathcal
B}(u_m+\xi),u'_m)_{L^2(D)}=0,
$$
and so, by using Cauchy inequality and Young inequality,
$$
|u'_m|_H^2+\frac12\frac{d}{dt}||u_m||^2_V \le
\frac12|u'_m|^2_H+\frac12\int_D |w|^2|\nabla u_m+\nabla\xi|^2,
$$
that is
$$ |u'_m|_H^2+\frac{d}{dt}||u_m||^2_V \le
\int_D |w|^2|\nabla u_m+\nabla\xi|^2,
$$
in particular
$$
\frac{d}{dt}\|u_m\|^2_V\le \int_D |w|^2|\nabla u_m+\nabla z|^2
\le 2\|w\|_{L^\infty(D)}^2\|u_m\|_V^2+2\|w\|_{L^\infty(D)}^2\|\xi\|_V^2.
$$
By Gronwall lemma we get
$$
\|u_m(t)\|_V^2
\le \|u_m^0\|^2_V{\rm e}^{2T\|w\|_{L^\infty(D_T)}^2}
+2\|w\|_{L^\infty(D_T)}^2\|\xi\|^2_{L^2(V)}{\rm e}^{2T\|w\|_{L^\infty(D_T)}^2},
$$
and, by integration by time, then we have
$$
\int_0^T|u'_m|_H^2\le
\|u_0\|_V^2+2\|w\|_{L^\infty(D_T)}^2\left(\|u_m\|_{L^\infty(V)}^2+\|\xi\|_{L^2(V)}^2\right).
$$
In conclusion we obtain that $u_m$ is bounded in $L^\infty(0,T;V)$
and $u'_m$ is bounded in $L^2(0,T;H)$. From the equation then we get
$$
Au_m=-u'_m-{\mathcal PB}(u_m+\xi)
$$
and so $Au_m\in L^2(0,T;H)$; by
the regularity theory for the Stokes operator, we obtain a bound
for $u_m$ in the space $L^2(0,T;H^2( D ))$.
\par
Then there exist a subsequence $(u_{m'})_{m'\in{\bf N}}$ of $(u_m)_{m\in{\bf N}}$ and
a function $u$ such that $u_{m'}$ converges weakly to $u$ in $L^2(0,T;V)$
and $L^2(0,T;H^2( D ))$ and converges weakly$^*$ in $L^\infty(0,T;V)$ and in
$L^\infty(0,T;H)$.
Moreover $u'_{m'}$ converges weakly to $u'$ in $L^2(0,T;H)$.

Taking the limit in the equation gives
$$
\frac{d}{dt}(u,v)_H+(u,v)_V+({\mathcal B}(u+\xi),v)_{L^2(D)}=0\qquad\mbox{for each }v\in V\\
$$
in the sense of distributions on $[0,T]$.

We can easily see that the solution is unique, that $u\in C([0,T];V)$ and
that
$$
u'+Au+{\mathcal P}{\mathcal B}(u+\xi)=0\qquad\text{in }H.
$$
This means that ${\mathcal P}\left(u'+Au+B(u+\xi)\right)=0$. Since 
$u'+Au+B(u+\xi)\in L^2(D_T)$, there exists a function $p$ such that 
$\nabla p(t)\in L^2(D)$ for a.e. t, and
$$
\partial_t u-\triangle u+(w\cdot\nabla)(u+\xi)+\nabla p=0
$$
whence $\nabla p\in L^2(D_T)$. Normalising $p$ by
imposing that $\int_D  p\,dx=0$, we obtain $p\in L^2(0,T;H^1(D))$.

Now we prove the energy equality. Let $\phi\in C^\infty_c(D_T )$
and $G=-(w\cdot\nabla)(u+\xi)$. Then
$$
\partial_t u-\triangle u+\nabla p=G;
$$
we mollify in ${\bf R}^4$ this equation in order to obtain smooth functions
$u_m$, $p_m$ and $G_m$ such that
$$
\begin{cases}
\partial_t u_m-\triangle u_m+\nabla p_m=G_m\\
\Div u_m=0
\end{cases}
$$
in a neighbourhood of $\supp\phi$ and such that
$$
\begin{aligned}
&u_m\to u&\qquad&\text{in }L^\infty(L^2(D)),\\
&\nabla u_m\to \nabla u&\qquad&\text{in }L^2,\\
&p_m\to p&\qquad&\text{in }L^2,\\
&G_m\to G&\qquad&\text{in }L^2;
\end{aligned}
$$
then we multiply by $u_m \phi$ and integrate by parts to have
$$
2\dint|\nabla u_m|^2\phi=
\dint|u_m|^2(\partial_t\phi+\triangle\phi)
      +2\dint p_m(u_m\cdot\nabla\phi)
      +2\dint (u_m\cdot G_m)\phi.
$$
As $m\to\infty$, we recover the energy equality, using the fact that
$G=-(w\cdot\nabla)(u+\xi)$ and by integration by parts.
\end{proof}

In the second step of the proof of Theorem \ref{thdetsws}, we obtain the
solution for the approximated nonlinear equation. We firstly define a 
regularisation procedure. Let $v_1,\ldots,v_m,\ldots$ be an orthonormal 
basis in $H$ of eigenfunctions of the operator $A$. For any $N\in{\bf N}$
and $v\in H$, we denote by $v^N$ the projection of $v$
on the span of $v_1,\ldots,v_N$. The following properties hold
\begin{enumerate}
\item $|v^N|_H\le |v|_H$;
\item $\|v^N\|_V\le \|v\|_V$;
\item $\|v-v^N\|_V\to 0$ if $v\in V$.
\end{enumerate}

Notice that, by virtue of assumption \eqref{Ad}, we know that
$$
z\in L^\infty(0,T;H)\cap L^2(0,T;V)\cap L^{8+\e}(0,T;L^4(D)).
$$
for some $\e>0$. Actually we know much more, namely that $z$ is bounded 
with values in $L^4(D)$, but, as we shall see, a weaker bound (like the
one given above) is sufficient.

\begin{lemma}\label{secondstep}
Assume \eqref{Ad}. Let $N\in\n$, then there exists a pair $(v_N,\pi_N)$,
with $v_N\in L^2(0,T;H^2(D))\cap C([0,T];H^1_0(D))$ and $\partial_t v_N\in L^2(D_T)$,
$\pi_N\in L^2(D_T)$ with $\int_D \pi_N=0$, that solves the 
following equation
$$
\partial_t v_N-\triangle v_N+\left[\left((v_N)^N+z^N\right)\cdot\nabla\right](v_N+z)+\nabla\pi_N=0
$$
with initial condition $v_N(0)=u_0^N$, and such that for any 
$\phi\in C^\infty_c(D_T )$, $\phi\ge 0$, the following energy
equality holds
\begin{align*}
2\intDt |\nabla v_N|^2\phi
&=      \intDt |v_N|^2(\partial_t\phi+\triangle\phi)
       +2\intDt  \pi_N(v_N\cdot\nabla\phi)\\
&\quad +\intDt (|v_N|^2+2v_N\cdot z)\left(((v_N)^N+z^N)\cdot\nabla\phi\right)\\
&\quad +2\intDt \phi z\cdot\left(((v_N)^N+z^N)\cdot\nabla\right)v_N
\end{align*}
\end{lemma}
\begin{proof}
Fix $N\in\n$ and let
$$
{\mathcal C}=\{w\in L^2(0,T;V)\cap L^\infty(0,T;H)\,|\,\|w\|_{\mathcal C}\le R_0\},
$$
where
$$
\|\cdot\|_{\mathcal C}=(\|\cdot\|_{L^2(V)}^2+\|\cdot\|_{L^\infty(H)}^2)^{1/2}
$$
and $R_0$ will be fixed later. Define a function ${\mathcal F}$ from 
${\mathcal C}$ to $L^\infty(0,T;H)\cap L^2(0,T;V)$ as follows:
if $w\in {\mathcal C}$, we take the regularisation $w^N$ as above and 
$u={\mathcal F}w$ is the solution of the problem
\begin{enumerate}
\item $\partial_t u-\triangle u+\left[(w^N+z^N)\cdot\nabla\right](u+z)+\nabla p=0$,
\item $\Div u=0$,
\item $u(0)=u_0^N$,
\item $u\in L^2(0,T;H^2( D ))\cap C([0,T];H^1_0( D ))$,
\item $\partial_t u$, $p\in L^2(D_T)$ \quad and \quad $\int_D p\,dx=0$,
\item for any $\phi\in C^\infty_c(D_T )$, $\phi\ge 0$,
\begin{align*}
2\intDt |\nabla u|^2\phi
&= \intDt |u|^2(\partial_t\phi+\triangle\phi)
       +2\intDt  p(u\cdot\nabla\phi)\\
&\quad +\intDt (|u|^2+2u\cdot z)\left((w^N+z^N)\cdot\nabla\phi\right)\\
&\quad +2\intDt \phi z\cdot\left((w^N+z^N)\cdot\nabla\right)u
\end{align*}
\end{enumerate}
The existence and uniqueness of this solution is guaranteed by the previous 
lemma, once we apply it with $w\to w^N+z^N$ and $\xi\to z$.

First we show that ${\mathcal F}$ maps ${\mathcal C}$ into itself. In 
order to show this, we shall only choose a suitable $R_0$. We have
\begin{align*}
\|{\mathcal F}w\|_{L^\infty(H)}^2+2\|{\mathcal F}w\|_{L^2(V)}^2
&\le |u_0^N|^2+2\intDt |z|\cdot|w^N+z^N|\cdot|\nabla{\mathcal F}w|\\
&\le |u_0|^2+\|{\mathcal F}w\|_{L^2(V)}^2+\intDt |z|^2|w^N+z^N|^2
\end{align*}
and so, since in finite-dimensional spaces all the norms are equivalent, we have
$$
\|w^N+z^N\|_{L^\infty}\le C_N\|w^N+z^N\|_{L^\infty(H)}\le C_N\|w+z\|_{L^\infty(H)}.
$$
Then
\begin{align*}
\|{\mathcal F}w\|_{L^\infty(H)}^2+2\|{\mathcal F}w\|_{L^2(V)}^2
&\le |u_0|_H^2+C_N^2T\|z\|_{L^\infty(H)}^2\|w+z\|^2_{L^\infty(H)}\\
&\le |u_0|_H^2+2C_N^2T\|z\|_{L^\infty(H)}^2\left(R_0^2+\|z\|^2_{L^\infty(H)}\right)\\
&\le R_0^2
\end{align*}
if we choose $R_0>|u_0|^2$ and $T$ small enough.

Then we show that ${\mathcal F}$ is a contraction. Let $w_1$, 
$w_2\in{\mathcal C}$ and set $w=w_1-w_2$, $v={\mathcal F}w_1-{\mathcal F}w_2$
and $p=p_1-p_2$, where $p_1$, $p_2$ are the corresponding pressures.
Then
$$
\partial_t v-\triangle v+\nabla p
=-(w^N\cdot\nabla)({\mathcal F}w_1+z)-\left[(w_2^N+z^N)\cdot\nabla\right]v,
$$
and so, using the fact that v(0)=0, we have
\begin{align*}
\|v\|^2_{L^\infty(H)}+2\|v\|^2_{L^2(V)}
&\le 2\intDt  |{\mathcal F}w_1+z|\cdot|w^N|\cdot|\nabla v|\\
&\le 2\|{\mathcal F}w_1+z\|_{L^\infty(H)}\cdot\|w^N\|_{L^\infty(D_T)}\int^T_0\|\nabla v\|_{L^2(D)}\\
&\le 2T^{\frac12}\|{\mathcal F}w_1+z\|_{L^\infty(H)}\cdot\|w^N\|_{L^\infty(D_T)}\|\nabla v\|_{L^2(D_T)}. 
\end{align*}
In particular
$$
\|v\|_{L^2(V)}
\le T^{1/2}\|{\mathcal F}w_1+z\|_{L^\infty(H)}\cdot\|w^N\|_{L^\infty(D_T)},
$$
and it follows that
\begin{align*}
\|v\|^2_{L^\infty(H)}+\|v\|^2_{L^2(V)}
&\le 2T\|{\mathcal F}w_1+z\|_{L^\infty(H)}^2\cdot\|w^N\|_{L^\infty(D_T)}^2\\
&\le 4TC_N^2\left(R_0^2+\|z\|^2_{L^\infty(H)}\right)\|w\|_{\mathcal C}^2.
\end{align*}
In conclusion
$$
\|{\mathcal F}w_1-{\mathcal F}w_2\|_{\mathcal C}^2\le 4TC_N^2\left(R_0^2+\|z\|^2_{L^\infty(H)}\right)\|w_1-w_2\|_{\mathcal C}^2
$$
and, if we choose the time interval small enough, the map ${\mathcal F}$ is a
contraction.
\end{proof}

Then the last step of the proof follows. We show that the sequence $(v_N,\pi_N)$
converges to a weak solution $(v,\pi)$ satisfying the properties of Definition
\ref{detsws}.

\begin{proof}[Proof of Theorem \ref{thdetsws}]
First we get an estimate of the solutions in the spaces $L^\infty(0,T;H)$
and $L^2(0,T;V)$. Indeed, multiply the equation by $v^N$ and
integrate by parts to get 
\begin{equation}\label{Ncei}
\frac12\frac{d}{dt}|v_N|_H^2+\|\nabla v_N\|_{L^2(D)}^2
= \int_D  z\cdot\left[\left((v_N)^N+z^N\right)\cdot\nabla\right]v_N,
\end{equation}
By using H\"older inequality and Young inequality we get
\begin{align*}
&\int_D  z\cdot\left[\left((v_N)^N+z^N\right)\cdot\nabla\right]v_N\le\\
&\qquad\le \int_D |z|\cdot|\nabla v_N|\cdot|(v_N)^N|
          +\int_D |z|\cdot|z^N|\cdot|\nabla v_N|\\
&\qquad\le \|\nabla v_N\|_{L^2(D)}\|(v_N)^N\|_{L^4(D)}\|z\|_{L^4(D)}
          +\|\nabla v_N\|_{L^2(D)}\|z\|_{L^4(D)}\|z^N\|_{L^4(D)}\\
&\qquad\le C^2\|\nabla v_N\|_{L^2(D)}^{7/4}|v_N|_{H}^{1/4}\|z\|_{L^4(D)}
          +C\|\nabla v_N\|_{L^2(D)}|z|_{H}^{1/4}\|z\|_{V}^{3/4}\\
&\qquad\le \frac12\|\nabla v_N\|_{L^2(D)}^2
          +C^2\|z\|^2_{L^4(D)}|z|_H^{1/2}\|z\|_V^{3/2}
          +\frac{7^7}{2^{10}}C^{16}\|z\|^8_{L^4(D)}|v_N|_H^2
\end{align*}
and so 
$$
\frac{d}{dt}|v_N|^2+\|\nabla v_N\|^2
\le 2C^2\|z\|_{L^4(D)}^2|z|^{1/2}\|z\|^{3/2}
    +\frac{7^7}{2^9}C^{16}\|z\|^8_{L^4(D)}|v_N|^2,
$$
since, for any suitable $\xi$, by virtue of Sobolev inequalities,
$$
\|\xi^N\|_{L^4(D)}
\le C|\xi^N|_H^{1/4}\|\xi^N\|_V^{3/4}
\le C|\xi|_H^{1/4}\|\xi\|_V^{3/4}.
$$
Then, by Gronwall lemma,
$$
\sup_{(0,T)}|v_N|_H^2\le C(T,z)+|u^N_0|^2_H\le |u_0|^2_H+C(T,z)
$$
and, integrating with respect to time,
$$
\int_0^T\|\nabla v_N\|^2_{L^2}\le |u_0|_H^2+C(T,z),
$$
where $C(T,z)$ is a constant which depends only on $T$ and on the function $z$. 

Then we give an estimate of the pressure term. By Theorem $15$ of 
\cite{Sol} we can deduce that $\nabla \pi_N$ are bounded in 
$L^{5/4}((\e,T);L^{5/4}(D))$ for every $\e>0$. Then
using the argument given in \cite{CKN} (page 781), we can conclude
that $\pi_N$ are bounded in $L^{5/4}((\e,T);L^{5/3}_{\text{loc}}(D))$,
provided that
$$
\int_D \pi_N\,dx=0
$$
at each time.

We can improve the regularity of $\pi_N$ using the general result of 
Sohr and Von Wahl \cite{SoVoW} or the simplified argument of Lin \cite{Lin},
to obtain that $\pi_N$ are bounded in $L^{5/3}_{\text{loc}}((0,T]\times D)$.

At last, using an argument similar to the one in Lemma 4.2 (Chapter III) 
of \cite{Te1}, we know that $v_N$ are bounded in $W^{1,2}(0,T;D(A^{-1}))$ 
and so, by virtue of Theorem 2.1 (Chapter III) of \cite{Te1}, 
$(v_N)_{N\in{\bf N}}$ is compact in $L^2(D_T )$.

We can deduce then that there exist a subsequence of $(v_N,\pi_N)_{N\in{\bf N}}$,
which we call again $(v_N,\pi_N)$, and functions $(v,\pi)$ such that
\begin{enumerate}
\item $v_N\to v$ weakly$^*$ in $L^\infty(0,T;L^2( D ))$,
\item $\nabla v_N\to\nabla v$ weakly in $L^2(D_T)$,
\item $v_N\to v$ strongly in $L^2(D_T)$,
\item $\pi_N\to \pi$ weakly in $L^{5/3}_{\rm loc}((0,T]\times D )$,
\item $v'_N$ is bounded in $L^2(0,T;D(A^{-1}))$.
\end{enumerate}

These convergence properties are sufficient to verify that the limit 
$v$ is a weak solution of Navier-Stokes system. Moreover the initial 
condition is satisfied, in fact $v_N$ are weakly continuous uniformly, 
by the bound of their derivatives, and so
$$
v(0)=\lim_Nv_N(0)=u_0.
$$
At last, thanks to the uniform bound on the time derivative, the limit 
is continuous as a function from $[0,T]$ to the space $H$ with the weak 
topology.

Now we prove the classical energy inequality. Integrate \eqref{Ncei} in 
time between $s$ and $t$, then in the limit as $N\to\infty$, the classical 
energy inequality for $v$ is obtained.

The last step of the proof is to prove that the limit $v$ verifies the 
local energy inequality. Since $v_N\to v$ in $L^2(D_T )$ and $v_N$ are 
bounded in $L^{10/3}(D_T )$ by Sobolev inequalities, then $v_N\to v$ 
in $L^q(D_T )$ for any $q\in [2,\frac{10}{3})$. By the properties of 
the regularisation, we can deduce that $(v_N+z)^N\to(v+z)$ in $L^2(D_T)$
and, in the same way as above, in $L^q(D_T)$.

Let $\varphi\in C^\infty_c(D_T)$. We know that
\begin{align*}
\dint|\nabla v_N|^2\varphi
&=      \dint|v_N|^2(\partial_t\varphi+\triangle\varphi)
       +\dint 2\pi_Nv_N\cdot\nabla\varphi\\
&\quad +\dint (|v_N|^2+2v_N\cdot z)\left((v_N+z)^N\cdot\nabla\varphi\right)\\
&\quad +2\dint\varphi z\cdot\left[(v_N+z)^N\cdot\nabla\right]v_N.
\end{align*}
By lower semi-continuity
$$
\dint|\nabla v|^2\varphi\le\liminf\dint|\nabla v_N|^2\varphi;
$$
moreover, since $v_N$ converges strongly in $L^q(D_T)$, with
$q\in[2,\frac{10}3)$, and $p_N$ converges weakly in $L^{5/3}_{\text{loc}}(D_T)$,
the first three terms converge. In order to show that the last term 
also converges, we use the fact that $z$ is bounded in $L^{8+\e}(0,T;L^4(D))$,
with $\e>0$ (this is the only step of the proof where we need this fact). Let
$$
q=\frac{6(8+\e)}{12+\e}\in(4,6),
\qquad p=\frac{4q}{3q-6}\in(2,\frac83)
$$
then we know that $(v_N+z)^N$ converges in $L^p(D_T)$, by the previous
considerations, and is bounded in $L^p(0,T;L^q(D))$ by the Sobolev
inequality. Thus by interpolation $(v_N+z)^N$ converges in the space 
$L^p(0,T;L^4(D))$ and this is sufficient to conclude since
$$
\frac12+\frac1p+\frac1{8+\e}=1.
$$
\end{proof}


\subsection{The proof of Theorem \ref{thmar}}\label{marex}
We can use now the results of the previous section to show the
existence of martingale suitable weak solutions. There are some
technical points in the proof of this theorem, mostly linked to
the fact that we deal with a pair $(v,\pi)$ of processes, where
we have no information on the tightness of the laws of the
approximating sequence of pressures. We solve the problem by means
of the following lemma, which is given in a generalised setting.

Let $(S,d_S)$ and $(T,d_T)$ be two complete separable metric spaces and 
consider the product metric space $S\times T$, endowed of the product metric.
Let $\pi:S\times T\to S$ be the canonical projection onto the first component,
that is $\pi(s,t)=s$ for $(s,t)\in S\times T$. Let a sequence of measure 
$\nu_n$ be given on $S\times T$ such that
$$
\mu_n=\pi\nu_n\rightharpoonup\mu,
$$
where $\mu$ is a measure on $S$. 

\begin{lemma}\label{skorohod}
There exist a probability space $(\Omega,\mathcal{F},\Pro)$, a sequence of 
random variables $(X_n,Y_n)$ on $S\times T$, the laws of which are 
$\nu_n$, and a random variable $X$ on $S$ having $\mu$ as its law, 
such that
$$
X_n\to X\qquad \Pro-a.s.
$$
\end{lemma}

\begin{proof}
From Theoreme 1, $\S\,6$, No. 1 of Bourbaki \cite{Bou}, since $T$ is a
complete separable metric space, it may be homeomorphic-ally embedded as
a $G_\delta$ subset (a countable intersection of open sets), and so as
a Borel set, of a compact metrizable space $\tilde T$. So measures 
$\nu_n$ can be extended to measures $\tilde\nu_n$ in the space $\tilde T$:
in this way the sequence $(\tilde\nu_n)_{n\in\n}$ is tight and by Prohorov
theorem there exists a subsequence, called again $(\tilde\nu_n)_{n\in\n}$,
converging weakly to a measure $\tilde\nu$ on $\tilde T$. Obviously 
$\tilde\pi\tilde\nu_n=\mu_n$ and $\tilde\pi\tilde\nu=\mu$, where 
$\tilde\pi$ is the projection of the space $S\times\tilde T$ onto the
first component. By Skorohod theorem there exist a probability space 
$(\Omega,\mathcal{F},\Pro)$, random variables $(\tilde X_n,\tilde Y_n)$ on
$S\times \tilde T$, with laws $\tilde\nu_n$, and $(\tilde X,\tilde Y)$,
with law $\tilde\mu$, such that
$$
(\tilde X_n,\tilde Y_n)\to (\tilde X,\tilde Y)\qquad\Pro-\text{a.s.}
$$
Since $\tilde\nu_n(S\times T)=1$, the restrictions to the space $S\times T$
of the previous random variables $(X_n,Y_n)$ (notice that $X_n=\tilde X_n$),
have $\nu_n$ as their laws. Moreover, $\tilde X$ has $\mu$ as its law and
$X_n\to\tilde X$, $\Pro$-a.\ s.
\end{proof}

\begin{remark}
In \cite{Rom} it is given an alternative proof of this fact, showing that
actually the claim is true for the whole sequence, not only for a subsequence.
Since we will use the lemma together with a compactness argument, we don't
really need the complete result.
\end{remark}

\begin{proof}[Proof of Theorem \ref{thmar}]
Let $z$ be the stochastic process which is solution of system \eqref{lNS}.
We know that for a.e. $\omega\in\Omega$ we have that
$W(\omega)$ is in the space $C^{1/2-\e}([0,T];D(A^\delta))$
for all $\e<\delta$, and so by \cite{Fla} we can deduce that 
$z(\omega)$ satisfies \eqref{zreg}. So we can apply path-wise the results 
of the previous section. For a.e. $\omega\in\Omega$ we use
Lemma \ref{secondstep} to get for each $N\in{\bf N}$ a 
pair $(v_N,\pi_N)$. The map 
$$
\omega\in\Omega\mapsto(v_N(\omega),\pi_N(\omega),z(\omega),W(\omega))
$$
with values in 
$$
L^2(0,T;H)\times L^{5/3}_{\rm loc}(D_T) \times C([0,T];H)\times C_0([0,T];H)
$$
is measurable for any $T>0$. The random variable $v_N(t)$ is measurable
for almost each $t\ge 0$ since each vector field is continuous with 
values in $H$. In fact, if $\eta_\e$ are mollifiers, then
$$
\eta_\e*v_N(t)\to v_N(t)\qquad\text{in }H
$$
and $\eta_\e*v_N$ are measurable. In the same way, thanks to uniqueness
of the solutions, we can use a smarter regularisation, namely
$$
v_N^0+\frac1\e\int_0^t\text{e}^{\frac{s-t}\e}(v_N(s)-v_N^0)\,ds
$$
to show that $v_N$ is a progressively measurable process.

So it is well defined a random variable 
$(v_N(\cdot,\omega),\pi_N(\cdot,\omega),z(\cdot,\omega),W(\cdot,\omega))$
such that $z$ satisfies \eqref{lNS}, $(v_N,\pi_N)$ is the solution of the
approximated problem for almost each $\omega\in\Omega$, and $W$ is a 
Brownian motion.

Let $\nu_N$ be the law of $(v_N,\pi_N,z,W)$ in
$$
{\mathcal E}=L^2(0,T;H)\times L_{\rm loc}^{5/3}(D_T)\times C([0,T];H)\times C_0([0,T];H)
$$
and let $\mu_N$ be the projection of $\nu_N$ in the variable $v$, that is
the law of $v_N$. We want to show that the family of measures $\mu_N$ is
tight, that is for each $\e>0$ there exists a compact set $K_\e$ in $\mathcal E$
such that
$$
\mu_N(K_\e)\ge 1-\e\qquad N\in{\bf N}.
$$
We take
$$
K_\e=\left\{v\,|\,\|v\|_{L^\infty(0,T;H)}^2+\|v\|_{L^2(0,T;V)}^2+\|v\|_{H^1(0,T;D(A^{-1}))}\le C_\e\right\}.
$$
The set $K_\e$ is compact in $L^2(0,T;H)$ and, moreover,
$$
{\Pro}\left[v\not\in K_\e\right]\le\frac1{C_\e}\Es(\|v_N\|_{L^\infty(0,T;H)}^2+\|v_N\|_{L^2(0,T;V)}^2+\|v_N\|_{H^1(0,T;D(A^{-1}))})
$$
and the right hand side is smaller than $\e$ if the above mean values are
uniformly bounded with respect to $N$ and $C_\e$ is chosen properly. To
see this, fix $N\in{\bf N}$ and let $u_N=v_N+z$. First we have (this can
be done as in Lemma 2.3 in \cite{Lin})
$$
\|\partial_t v_N\|_{L^2(0,T;D(A^{-1}))}
\le \|v_N\|_{L^2(0,T;V)}+\|u_N\|_{L^2(0,T;H)}\|u_N\|_{L^2(0,T;V)},
$$
moreover ${\Es}\|z\|_{L^2(0,T;V)}^2+{\Es}\|z\|^2_{L^\infty(0,T;H)}$ is 
finite and so the only thing we need to show is that
$$
{\Es}\|u_N\|_{L^\infty(0,T;H)}^2+{\Es}\|u_N\|_{L^2(0,T;V)}^2
$$
is bounded uniformly in $N$. We know that
$$
u_N(t)=v_0^N+\int_0^t(A u_N+((u_N)^N\cdot\nabla)u_N+f^N)\,ds+W_t
$$
and so (see Pardoux \cite{Par}, Th\'eor\`eme 3.1)
\begin{multline}\label{pardu}
|u_N(t)|^2_H+2\int_0^t\|u_N\|_V^2\,ds=\\
= |u_0^N|^2_H+2\int^t_0\langle f^N,u_N\rangle_H\,ds
 +2\int_0^t\langle u_N(s),dW_s\rangle_H+\sigma t.
\end{multline}
Notice that, if $\tau_R=\inf\{t>0\,|\,|u_N(t)|>R\}$, then
\begin{equation}\label{martin}
\int_0^t\langle u_N(s),dW_s\rangle_H
\end{equation}
is a local martingale with respect to the stopping time $\tau_R$,
and so, taking the expectation of \eqref{pardu} at time $t\wedge\tau_R$, 
\begin{multline*}
\Es|u_N(t\wedge \tau_R)|^2_H+2{\Es}\int_0^{t\wedge\tau_R}\|u_N\|_V^2\,ds=\\
= \Es|u_0^N|^2_H
 +\Es\int^{t\wedge\tau_R}_0|f^N|_H^2\,ds
 +{\Es}\int^{t\wedge\tau_R}_0|u_N|^2_H\,ds
 +\sigma t\wedge\tau_R,
\end{multline*}
Let $\varphi(t)={\Es}|u_N(t\wedge\tau_R)|^2_H$, then we have
$$
\varphi(t)\le\varphi(0)+\int^t_0|f|_H^2\,ds+\int_0^t\varphi(s)\,ds+\sigma t
$$
and by Gronwall's lemma we can deduce that $\varphi(t)$ is bounded by a
constant independent of $R$. So, as $R\uparrow\infty$, we can deduce that
$\Es|u^N(t)|^2_H\le C(T)$, and then that \eqref{martin} is a martingale.

So, by taking the expectation in \eqref{pardu}, we obtain first that
$$
\Es|u_N(t)|^2_H+\Es\int_s^t\|u_N\|_V^2\,dr
\le \Es|u_N(s)|^2_H 
    +\sigma (t-s)
    +\Es\int^t_s\|f^N\|_{V'}\,ds,
$$
and then, using the Burkholder-Davis-Gundy inequality, that
\begin{multline*}
\Es\bigl[\sup_{s\le t}|u_N(s)|^2_H\bigr]+\Es\int_0^t\|u_N\|_V^2\,ds\le\\
\le \Es|u_0^N|^2_H
   +\int^t_0\|f^N\|_{V'}^2\,ds
   +\sigma t
   +2\sigma C_1{\Es}\bigl[\int_0^t|u_N|^2_H\,ds\bigr]^{\frac12}.
\end{multline*}
We can conclude that
$$
\Es\bigl[\sup_{s\le t}|u_N(s)|^2\bigr]+\Es\int_0^t\!\!\|u_N\|^2\,ds
\le 2\Es|u_0^N|^2
    +2\int^t_0\!\!\|f^N\|_{V'}^2\,ds
    +2\sigma (1+\sigma C_1^2)\,t
$$
and the claim is proved.

By Lemma \ref{skorohod}, there exist a probability
space $(\tilde\Omega,\tilde{\mathcal F},\Pro)$ and random variables
$\tilde U_{N_k}=(\tilde v_{N_k},\tilde\pi_{N_k},\tilde z,\tilde W)\in{\mathcal E}$
such that the law of each $\tilde U_{N_k}$ is $\nu_{N_k}$ and
$$
\tilde v_{N_k}\to \tilde v\qquad\tilde\Pro-a.s.,
$$
where $\tilde v$ is a random variable whose law is $\mu$.

It is easy to check that $\tilde W$ is a Wiener process which keeps the 
same regularity properties of $W$. Notice that
$$
\Pro\left[v_N\in L^2(0,T;D(A))\cap C([0,T];V)\right]=1,
$$
and so the same holds true for the new random variables $\tilde v_{N_k}$.
In the same way we can deduce that $\pi_{N_k}\in L^2(D_T)$ and so on. 
Now we need to show that the $\tilde v_{N_k}$ satisfy the equations and 
the energy inequalities. We give a proof, using a trick of Bensoussan 
\cite{Ben}, for the local energy inequality (actually it is an equality 
for the $v_N$). Given $\phi\in C_c^\infty(D_T)$, define the random variable 
$X^N:(\Omega,{\mathcal F})\to ({\bf R},{\mathcal B}({\bf R}))$ as
\begin{align*}
X^N
&=      \Big|2\intDt|\nabla v_N|^2\varphi
       -\intDt|v_N|^2(\partial_t\varphi+\triangle \varphi)\\
&\quad -\intDt(|v_N|^2+2v_N\cdot z)\left((v_N+z)^N\cdot\nabla\varphi\right)\\
&\quad -2\intDt\varphi z\cdot((v_N+z)^N\cdot\nabla)v_N
       -\intDt 2\pi_N v_N\cdot\nabla\varphi\Big|,
\end{align*}
and let $\tilde X^{N_k}$ be the analogue of $X^N$ for the $\tilde U_{N_k}$.

We know that $X^N=0$, $\Pro$-a.\ s., and so
$$
{\Es}\frac{X^N}{1+X^N}=0.
$$
Notice that 
$$
\frac{X^N}{1+X^N}=\Phi(u_N),
$$
where $\Phi$ is a deterministic bounded continuous function on the 
subspace of $\mathcal E$ where the $\nu_N$ are concentrated (remember 
that the $U_N$ are far more regular than the elements of ${\mathcal E}$)
and so
$$
\tilde{\Es}\frac{\tilde X^{N_k}}{1+\tilde X^{N_k}}
=\tilde{\Es}\Phi(\tilde U_{N_k})
=\int_{\mathcal E}\Phi(u)\,\nu_{N_k}(du)
={\Es}\Phi(U_{N_k})
={\Es}\frac{X^{N_k}}{1+X^{N_k}}
=0;
$$
this means
$$
\tilde X_{N_k}=0\qquad \tilde{\Pro}-a.s.
$$
If we do this for a dense set of functions in $C_c^\infty(D_T)$, we can 
conclude that there exists a set $\tilde\Omega_0\subset\tilde\Omega$ of 
full measure such that the local energy inequality holds for each 
$\omega\in\tilde\Omega_0$ and $\phi\in C_c^\infty(D_T)$.

From now on, since the two sequences enjoy the same properties, we will 
omit the tilde.

The last step of the proof is to show that the limit process is a 
martingale solution. We need to find the limit of the sequence of the 
pressures in such a way that the equations and the energy inequalities 
are satisfied. First we observe that the $v_{N_k}$ solve the equation
$$
\partial_t v_{N_k}+Av_{N_k}+B((v_{N_k})^{N_k}+z,v_{N_k}+z)=0
$$
$\Pro$-a.s. in $V'$, so in the limit
$$
\partial_t v+Av+B(v+z,v+z)=0,
\qquad\text{$\Pro$-a.s in $V'$}.
$$
Thus there exists a distribution $\pi$ such that \eqref{mNS} holds true.
Normalise $\pi$ in such a way that
$$
\int_D\pi(t)\,dx=0\qquad\text{a.e. }t.
$$

The set of $\omega\in\Omega$ such that $W(\omega)$, and so $z(\omega)$
and all $v_{N_k}(\omega)$, has the suitable regularity we need, such 
that $(v_{N_k}(\omega),\pi_{N_k}(\omega))$ satisfy the modified 
Navier-Stokes equations and such that $v_{N_k}(\omega)\to v(\omega)$, 
has probability one. Take an $\omega\in\Omega$ in this way. Then there 
exists a subsequence of $\pi_{N_k}(\omega)$ which converges weakly in 
$L^{5/3}_{\rm loc}(D_T)$. Taking the limit in the equations, we observe 
that the equations are satisfied both by $\pi(\omega)$ and by the limit 
of $\pi_{N_k}(\omega)$. This means that the two are equal (they have 
both zero mean in $D$) and
$$
\pi_{N_k}(\omega)\to\pi(\omega)\qquad\text{\rm weakly in }L^{5/3}_{\rm loc}(D_T).
$$
We need only to verify that $(v(\omega),\pi(\omega))$ satisfies the local 
energy inequality. This can be done as in the third step of the proof of 
Theorem \ref{thdetsws}, since $v_{N_k}(\omega)$ converges to $v(\omega)$ 
strongly in $L^2(D_T)$, weakly in $L^2(0,T;H^1_0(D))$ and weakly${}^*$ in 
$L^\infty(0,T;L^2(D))$, while $\pi_{N_k}(\omega)$ converges to $\pi(\omega)$ 
weakly in $L^{5/3}_{\rm loc}(D_T)$.

Finally we set $u=v+z$ and $P=\pi+Q$ and we can conclude that $(u,P)$ is a 
martingale suitable weak solution in the sense of Definition \ref{marsws}.
\end{proof}


\subsection{The proof of Theorem \ref{thstat}}
In this last section we prove the existence of stationary solutions.
The proof is given using the classical Krylov-Bogoliubov method, where
the initial measure is given by the law of a martingale solution.

In order to show the existence of time-stationary measures, we need the
following compactness lemma.

\begin{lemma}\label{compatto}
Let $(T_N)$ be a sequence of positive real numbers such that 
$T_N\uparrow\infty$, let $k(T_N)$ be an increasing sequence of positive
constants and let $\beta$, $s$, $p>0$ be such that $\beta>0$, $s<\frac12$ 
and $s\,p<1$. Then the set $K$ of all $(u,W)\in{\mathcal S}$ such that
$$
\|u\|^2_{L^\infty(0,T_N;H)}+\|u\|^2_{L^2(0,T_N;V)}+||W||_{W^{s,p}(0,T_N;D(A^\beta))}^p\le k(T_N)
$$
for each $N\in\n$ is compact in ${\mathcal S}$.
\end{lemma}
\begin{proof}
Since the $u\in K$ are bounded in $L^\infty(0,T;H)$ and in  $L^2(0,T;V)$ 
and they satisfy equation \eqref{NS} in distributions, it follows that
they are bounded in $H^1(0,T;D(A^{-1})$ (see Temam \cite{Te1}).

Moreover the immersion of the space $W^{s,p}(0,T;D(A^\beta))$ in $C([0,T];H)$
is compact. In conclusion, for any given $T_N$, we need only to show that,
if $(u_n,W_n)\in K$ and
$$
(u_n,W_n)\to(u,W)\qquad\text{in }L^2(0,T_N;H)\times C_0([0,T_N];H)
$$
then $(u,W)\in{\mathcal S}$, that is $(u,W)$ is a suitable weak solution
in $[0,T_N]$.

Let $z_n$ be the solution of the Stokes equation \eqref{lNS} with 
$\partial_t W_n$ as a forcing term and let 
$v_n=u_n-z_n$ and $\pi_n=P_n-Q_n$. By well known results on the 
Stokes equation (see \cite{Fla}) we know that $z_n$ are bounded 
in $L^\infty(0,T;H)$, $L^2(0,T;V)$ and $L^\infty(0,T;L^4(D))$ and 
moreover $z_n\to z$ in $L^2(0,T;H)$. Then $v_n$ are bounded in 
$L^\infty(0,T;H)\cap L^2(0,T;V)$ and so we can proceed as in the 
third step of the proof of Theorem \ref{thdetsws} to get all the 
convergence properties we need to take the limit in the equations 
and in the local energy inequality.
\end{proof}

\begin{proof}[Proof of Theorem \ref{thstat}]
We use the Krylov-Bogoliubov procedure for the semigroup $\tau_t$ in 
${\mathcal S}$. Let $u_0\in H$. In the previous section we have shown 
the existence of at least one martingale suitable weak solution 
$\overline u$ of Navier-Stokes system driven by a Brownian motion 
$W$ and with initial condition $u_0$. Let $\nu_0\in M_1({\mathcal S})$ 
be the law of the stochastic process $(\overline u,W)$ with values in 
${\mathcal S}$. Let $\nu_t=\tau_t\nu_0$ and set
$\mu_t=\frac1t\int_0^t\nu_s\,ds$.
Notice that the $W$-component of $\nu_t$ and $\mu_t$ is always the 
Wiener measure given by the Brownian motion $W$, due to the stationarity 
of this process.

Suppose that for each $\e>0$ there exists a compact set $K_\e$ in 
${\mathcal S}$ such that
\begin{equation}\label{tight}
\nu_t(K_\e)\ge 1-\e\qquad\mbox{for all }t\ge0
\end{equation}
(this claim will be proved in the sequel of the proof), so that
$\mu_t(K_\e)\ge 1-\e$, 
and the family of measures $(\mu_t)_{t\ge0}$ is tight. By means of 
Prohorov theorem there is a subsequence $(\mu_{t_n})_{n\in{\bf N}}$ 
which converges weakly to some $\mu\in M_1({\mathcal S})$. The measure 
$\mu$ is time stationary, in fact if $t\ge 0$ and $\phi\in C_b({\mathcal S})$,
\begin{align*}
(\tau_t\mu)(\phi) 
&=\mu(\tau_t\phi)\\
&=\lim_{n\to\infty}\mu_{t_n}(\tau_t\phi)\qquad(\text{since }\tau_t\phi\in C_b(W))\\
&=\lim_{n\to\infty}\frac1{t_n}\int_0^{t_n}\nu_0(\tau_{s+t}\phi)\,ds\\
&=\lim_{n\to\infty}\frac1{t_n}\int_t^{t+t_n}\nu_0(\tau_r\phi)\,dr\\
&=\lim_{n\to\infty}\mu_{t_n}(\phi)
   +\lim_{n\to\infty}\frac1{t_n}\left(\int_{t_n}^{t+t_n}\nu_0(\tau_r\phi)\,dr
                                      -\int_0^t\nu_0(\tau_r\phi)\,dr\right)\\
&=\mu(\phi)\qquad(\nu_0(\tau_r\phi)\text{ is bounded in }r).
\end{align*}
Then we show the claim in \eqref{tight}. In order to show that $\nu_t$ is 
a tight family of measures, we need only to show that 
$$
\Es\bigl[\sup_{(0,T)}|\tau_t\overline u|^2_+\int_0^T\|\tau_t \overline u\|^2\,dt+||\tau_t W||_{W^{s,p}}^p\bigr]
\le C(T)\qquad\text{uniformly in }t\ge 0,
$$
in fact if we take $k(T_n)>2^{n+1}\e C(T_n)$ and $K_\e$ as in Lemma 
\ref{compatto}, we have
\begin{align*}
\nu_t[K_\e^c]
& = \Pro[\tau_t(\overline u,W)\not\in K_\e]\\
& = \Pro\bigl[\cup_{n\in{\bf N}}\{\sup_{(0,T_n)}|\tau_t\overline u|^2_+\int_0^{T_n}\|\tau_t \overline u\|^2\,dt+||\tau_t W||_{W^{s,p}}^p>k(T_n)\}\bigr]\\
&\le\sum_{n=0}^\infty\Pro\bigl[\sup_{(0,T_n)}|\tau_t\overline u|^2_+\int_0^{T_n}\|\tau_t \overline u\|^2\,dt+||\tau_t W||_{W^{s,p}}^p>k(T_n)\bigr]\\
&\le\sum_{n=0}^\infty\frac1{k(T_n)}\Es\bigl[\sup_{(0,T_n)}|\tau_t\overline u|^2_+\int_0^{T_n}\|\tau_t \overline u\|^2\,dt+||\tau_t W||_{W^{s,p}}^p\bigr]\\
&\le\sum_{n=0}^\infty\frac{C(T_n)}{k(T_n)}<\e.
\end{align*}
To prove the claim, we use \eqref{energyinequality} and Poincar\'e
inequality to get
$$
{\Es}|\overline u(t)|^2_H+{\lambda}{\Es}\int_s^t|\overline u(r)|^2_H\,dr\le {\Es}|\overline u(s)|^2_H+(\sigma+\|f\|_{V'}^2)(t-s),
$$
then using Lemma \ref{gronwall} it follows that
$$
{\Es}|\overline u(t)|^2\le \sup_{s\in(0,1)}{\Es}|\overline u(s)|^2+\frac{\sigma+\|f\|_{v'}}{{\lambda}^2}
\qquad\text{for almost every $t$}.
$$
Finally by \eqref{supenergy} we obtain
\begin{align}\label{uniforme}
\Es|\tau_t \overline u|^2_{L^\infty(0,T;H)\cap L^2(0,T;V)}
& =      \Es\bigl[\sup_{(t,t+T)}|\overline u(s)|^2_H+\int_t^{t+T}\|\overline u(s)\|_V^2\,ds\bigr]\notag\\
&\le    2\Es|\overline u(t)|_H^2+2\left(\|f\|_{V'}^2+\sigma(1+\sigma C_1^2)\right)T\notag\\
&\le    2\sup_{s\in(0,1)}{\Es}|\overline u(s)|^2+\frac{2\sigma+\|f\|_{V'}}{{\lambda}^2}\\
&\quad +2\left(\|f\|_{V'}^2+\sigma(1+\sigma C_1^2)\right)T.\notag
\end{align}

The estimate on the Brownian motion is classical:
\begin{align*}
\Es\|\tau_tW\|_{W^{s,p}}^p
&\le    \int^T_0{\Es}\|\tau_tW(r)\|_{D(A^\beta)}^p\,dr\\
&\quad +\int_0^T\int_0^T\frac{{\Es}\|\tau_tW(r_1)-\tau_tW(r_2)\|_{D(A^\beta)}^p}{|r_1-r_2|^{1+sp}}\,dr_1\,dr_2\\
&\le    C_p\int_0^Tr^{p/2}\,dr+C_p\int^T_0\int^T_0|r_1-r_2|^{(\frac12-s)p-1}\,dr_1\,dr_2\\
&\le    C_{p,T},
\end{align*}
since $s<\frac12$.

We want to show now that the stationary measure has finite mean 
dissipation rate. We consider
$$
\phi(u)=\int^T_0\|u(t)\|^2\,dt,
$$
this is a lower semi-continuous function on ${\mathcal S}$, then 
there exists an increasing sequence of functions 
$\phi_N\in C_b({\mathcal S})$ such that $\phi_N\uparrow\phi$. 
From the monotone convergence theorem $\langle\phi_N,\mu\rangle$ 
converges to $\langle\phi,\mu\rangle$, even if the last term is 
not finite. So it is sufficient to show that $\langle\phi_N,\mu\rangle$ 
is bounded independently from $N$. Now
$$
\langle\phi_N,\mu\rangle=\lim_{n\to\infty}\frac1{t_n}\int^{t_n}_0\langle\phi_N,\nu_s\rangle\,ds
$$
and by \eqref{uniforme},
$$
\langle\phi_N,\nu_s\rangle
=   \Es\phi_n(\tau_s(\overline u))
\le \Es\phi(\tau_s(\overline u))\\
\le \Es\|\tau_s(\overline u)\|_{L^2(0,T;V)}
\le  C_T.
$$

Finally we show \eqref{lineare}. Let
$$
\Theta(t)=\int_{\mathcal S}\left[\int_0^t\|u\|^2\right]\mu(du),
$$
then by the invariance of $\mu$,
\begin{align*}
\Theta(t)-\Theta(s)
&= \int_{\mathcal S}\int_0^{t-s}\|u(r+s)\|^2\,dr\,\mu(du)\\
&= \int_{\mathcal S}\int_0^{t-s}\|u(r)\|^2\,dr\,\mu(du)\\
&= \Theta(t-s).
\end{align*}
Since $\Theta$ is non decreasing, then $\Theta(t)=Ct$.  
\end{proof}

\begin{remark}
When the dynamic is well defined, one can be interested in studying 
other mathematical objects, which can give some asymptotic information 
on the solutions. For example the dynamic for the linear Stokes equations 
is well defined and one can study the invariant measures of this equation. 
Then it can be easily seen that any time-stationary solution in the path 
space, {\sl frozen} at an arbitrary time, is an invariant measure.
In fact let
$$
p_T:z\mapsto z(T):C([0,\infty);H)\to H,
$$
such mapping is continuous. Let $\mu^z$ be the the time-invariant measure 
which can be built for the Stokes equation.
\end{remark}

\begin{proposition}
The image measure of $\mu^z$ through $p_T$ is an invariant measure for 
the Stokes equation \eqref{lNS}.
\end{proposition}
\begin{proof}
By the proof of the previous theorem
$$
\mu^z=\lim_{n\to\infty}\frac1{t_n}\int_0^{t_n}\tau_s\nu_0^z\,ds,
$$
where $\nu_0^z$ is the law of $\omega\to z(\cdot,\omega)$ in 
$C([0,\infty];H)$. Then
$$
p_T\mu^z=\lim_{n\to\infty}\frac1{t_n}\int_0^{t_n}p_T\tau_s\nu_0^z\,ds
        =\lim_{n\to\infty}\frac1{t_n}\int_0^{t_n}p_{T+s}\nu_0^z\,ds,
$$
since $p_T\circ\tau_s=p_{T+s}$. Now, since $\nu_0^z$ is the law of 
$z(\cdot,0)$, then $p_{T+s}\nu_0^z$ is the law of $z(T+s,0)$, that 
is the law of $z(s,z(T,0))$. In conclusion
$$
p_T\mu^z=\lim_{n\to\infty}\frac1{t_n}\int_0^{t_n}{\mathcal L}z(s,z(T,0))\,ds
$$
and, by Proposition $11.3$ of \cite{DaPZa}, $p_T\mu^z$ is an invariant 
measure.
\end{proof}

Finally we prove the easy exotic Gronwall lemma we used in the proof 
of the previous theorem.

\begin{lemma}\label{gronwall}
Suppose the function $v:[0,+\infty)\to{\bf R}$ satisfies
$$
v(t)\le v(s)-{\lambda}\int_s^tv(r)\,dr+C(t-s)
$$
for all $t\ge 0$ and almost all $s\le t$. Then
$$
v(t)\le \sup_{s\in(0,1)}v(s)+\frac C{\lambda}
$$
for almost all $t\ge 0$.
\end{lemma}
\begin{proof}
Let $t>0$ and let ${\mathcal N}$ be the set of Lebesgue measure zero 
for which the inequality does not hold. Set $u(s)=-v(t-s)$ for $s\in[0,t]$.
It is easy to see that for each $s\in[0,t]$ such that $t-s\not\in{\mathcal N}$,
we have
$$
u(s)\le u(0)+{\lambda}\int_0^s u(r)\,dr+Cs
$$
and, by Gronwall lemma
$$
u(s)\le u(0)\text{e}^{{\lambda} s}+\frac{C}{{\lambda}}(\text{e}^{{\lambda} s}-1).
$$
This means
\begin{equation*}
v(t)
\le v(t-s)\text{e}^{-{\lambda} s}+\frac{C}{{\lambda}}(1-\text{e}^{-{\lambda} s})
\le v(t-s)+\frac{C}{{\lambda}},
\end{equation*}
and then we can conclude that
$$
v(t)\le \sup_{s\in(0,1)}v(s)+\frac C{\lambda}.
$$
\end{proof}

\bibliographystyle{amsplain}

\begin{thebibliography}{99}
\bibitem{Ben}\textsc{A.~Bensoussan}, \textsl{ Stochastic Navier-Stokes
  equations}, Acta Appl. Math. {\bf 38} (1995), 267-304.
\bibitem{Bei}\textsc{H.~Beirao~da~Veiga}, \textsl{ On the construction
 of suitable weak solutions to the Navier Stokes equations via a general
 approximation theorem}, J. Math. Pures Appl., IX. Ser. 64 (1985),
 321-334.
\bibitem{Bou}\textsc{N.~Bourbaki}, \textsl{ Topologie g\'en\'erale},
  \'El\'ements de Math\'ematique, Hermann, Paris (1958).
\bibitem{CKN}\textsc{L.~Caffarelli, R.~Kohn, L.~Nirenberg}, \textsl{ Partial
  regularity of suitable weak solutions of the Navier-Stokes equations},
  Comm. Pure Appl. Math. XXXV (1982), 771-831.
\bibitem{DaPZa}\textsc{G.~Da~Prato, J.~Zabczyk}, \textsl{ Stochastic
  equations in infinite dimension}, Cambridge Univ. Press, Cambridge
  1992.
\bibitem{Fla}\textsc{F.~Flandoli}, \textsl{ Stochastic differential
  equations in fluid dynamics}, Rendiconti del Seminario
  Fisico-Matematico di Milano (1996). 
\bibitem{FlRo1}\textsc{F.~Flandoli, M.~Romito}, \textsl{ Statistically
  stationary solutions to the 3D Navier-Stokes equations do not show
  singularities}, Elec. J. Prob. \textbf{6} (2001).
\bibitem{FlRo2}\textsc{F.~Flandoli, M.~Romito}, \textsl{ Partial regularity
  for the stochastic Navier-Stokes equations},
  Trans. Amer. Math. Soc. \textbf{354}, no. 6 (2002), 2207--2241.
\bibitem{FlSc}\textsc{F. Flandoli, B. Schmalfuss}, \textsl{ Weak solutions
  and attractors for the 3D Navier-Stokes equations with non-regular
  force}, J. Dynam. Diff. Eq. {\bf 11}, Nr. 2 (1999), 355-398.
\bibitem{IkWa}\textsc{N. Ikeda, S. Watanabe}, \textsl{ Stochastic
  differential equations and diffusion processes}, North Holland
  Mathematical Library {\bf 24}, North Holland/Kodansha (1989).
\bibitem{Lem}\textsc{P. G. Lemarie-Rieusset}, \textsl{ Solutions faibles
  d'energie infinie pour les equations de Navier-Stokes dans ${\bf
    R}^3$}, C. R. Acad. Sci., Paris Ser. I, Math. {\bf 328}, No. 12
  (1999), 1133-1138.
\bibitem{Lin}\textsc{F. Lin}, \textsl{ A new proof of the
  Caffarelli-Kohn-Nirenberg theorem}, Comm. Pure Appl. Math. LI (1998),
  241-257.
\bibitem{Lio}\textsc{P. L. Lions}, \textsl{ Mathematical topics in fluid
  dynamic}, Vol. 1, Clarendon Press, Oxford (1996).
\bibitem{Par}\textsc{E. Pardoux}, \textsl{ Equations aux derivees partielles
  stochastiques nonlineaires monotones. Etude de solutions fortes de
  type Ito}, These, Universit\'e Paris Sud, Novembre 1975.
\bibitem{Rom}\textsc{M. Romito}, \textsl{ Partial regularity theory for a
  stochastic Navier-Stokes system}, Thesis, Pisa (2000).
\bibitem{Sc1}\textsc{V. Scheffer}, \textsl{ Partial regularity of solutions
  to the Navier-Stokes equations}, Pacific J. Math. {\bf 66} (1976),
  532-552.
\bibitem{Sc2}\textsc{V. Scheffer}, \textsl{ Hausdorff measure and the
  Navier-Stokes equations}, Comm. Math. Phys. {\bf 61} (1978), 41-68.
\bibitem{Sc3}\textsc{V. Scheffer}, \textsl{ The Navier-Stokes equation on a
  bounded domain}, Comm. Math. Phys. {\bf 73} (1980), 1-42.
\bibitem{Sel}\textsc{G. Sell}, \textsl{ Global attractor for the 3D
  Navier-Stokes equations}, J. Dynam. Diff. Eq. {\bf 8} (1), 1996.
\bibitem{Sol}\textsc{V. A. Solonnikov}, \textsl{ Estimates of the solutions
  of a nonstationary linearized system of Navier-Stokes equations},
  Amer. Math. Soc. Translations, Ser. 2, Vol. 75, 1-117.
\bibitem{SoVoW}\textsc{H. Sohr, W. von Wahl}, \textsl{ On the regularity of
  the pressure of weak solutions of Navier-Stokes equations},
  Arch. Math. {\bf 46}, Basel (1986), 28-439.
\bibitem{Te1}\textsc{R. Temam}, \textsl{ The Navier-Stokes Equations}, North
  Holland, 1977.
\bibitem{Te2}\textsc{R. Temam}, \textsl{ Navier-Stokes Equations and
  Nonlinear Functional Analysis}, SIAM, Philadelphia, 1983.
\bibitem{Te3}\textsc{R. Temam}, \textsl{ Infinite Dimensional Dynamical
  Systems in Mechanics and Physics}, Springer-Verlag, New York, 1988.
\end{thebibliography}
 
\end{document}